\documentstyle[twoside]{article}
\title{\bf Certain results on Euler-type integrals and their applications\\}
\author{\sc 
S.Jabee$^{a}$, M. Shadab$^{a}$\footnote{Corresponding author: M. Shadab} \ and R.B. Paris$^{b}$\footnote{E-mail: {\tt saimajabee007@gmail.com} (S. Jabee), {\tt shadabmohd786@gmail.com} (M. Shadab), {\tt r.paris@abertay.ac.uk} (R.B. Paris)}\\
\\
\\
$^{a}${\em Department of Applied Sciences and Humanities,}\\
{\em Faculty of Engineering and Technology,}\\
{\em Jamia Millia Islamia University, New Delhi 110025, India}\\
$^{b}${\em Division of Computing and Mathematics, Abertay University,}\\
{\em Dundee DD1 1HG, UK.}}

\begin{document}
\newcommand{\bee}{\begin{equation}}
\newcommand{\ee}{\end{equation}}
\def\f#1#2{\mbox{${\textstyle \frac{#1}{#2}}$}}
\def\dfrac#1#2{\displaystyle{\frac{#1}{#2}}}
\newcommand{\fr}{\frac{1}{2}}
\newcommand{\fs}{\f{1}{2}}
\newcommand{\g}{\Gamma}
\newcommand{\br}{\biggr}
\newcommand{\bl}{\biggl}
\newcommand{\ra}{\rightarrow}
\renewcommand{\topfraction}{0.9}
\renewcommand{\bottomfraction}{0.9}
\renewcommand{\textfraction}{0.05}
\date{}
\maketitle
\pagestyle{myheadings}
\markboth{\hfill {\it S. Jabee, M. Shadab and R.B. Paris} \hfill}
{\hfill {\it Certain results on Euler-type integrals } \hfill}
\begin{abstract} 
This paper deals with the evaluation of some definite Euler-type integrals in terms of the Wright hypergeometric function. We obtain a theorem on the Wright hypergeometric function and then use this theorem to evaluate some definite integrals. Further, we derive some results as applications of these evaluations. Multi-variable cases of the derived results of this paper are also briefly discussed.
\vspace{0.4cm}

\vspace{.2cm}


\end{abstract}

\vspace{0.2cm}

\noindent $\,$\hrulefill $\,$

\vspace{0.2cm}

\begin{center}
{\bf 1. \  Introduction}
\end{center}
\setcounter{section}{1}
\setcounter{equation}{0}
\renewcommand{\theequation}{\arabic{section}.\arabic{equation}}
The generalization of the generalized hypergeometric series ${}_pF_q$ is due to Fox \cite{Fox} and Wright \cite{Wright1,Wright2} who studied the asymptotic expansion of the so-called Wright hypergeometric function defined by (see \cite[p. 21]{Srivastava2})
\begin{equation}\label{pPsiq}
 {}_p\Psi_q \left[ \begin{array}{r}(\alpha_1,\,A_1),\,\ldots,\,(\alpha_p,\,A_p);\\
                      (\beta_1,\,B_1),\,\ldots,\,(\beta_q,\,B_q); \end{array}z \right] 
=\sum_{k=0}^\infty \,\frac{\prod_{j=1}^p\,\Gamma( \alpha_j + A_j\,k)}{\prod_{j=1}^q\, \Gamma ( \beta_j + B_j\,k)}\, \frac{z^k}{k!}.
    \end{equation}
The coefficients $A_j \in {\bf R}^+$ $(j=1,\ldots, p)$ and $B_j \in {\bf R}^+$ $(j=1,\ldots, q)$ are such that
\begin{equation}\label{Condi}
1 + \sum_{j=1}^q\, B_j - \sum_{j=1}^p\, A_j \geq 0.
    \end{equation}
Here and in the following, let ${\bf C}$, ${\bf R}$, ${\bf R}^+$,    ${\bf Z}$ and ${\bf N}$ be
the sets of complex numbers, real numbers, positive real numbers,  integers  and positive integers, respectively,
 and let
 $${\bf R}^+_0:={\bf R}^+ \cup \{0\}, \quad    {\bf N}_0:={\bf N} \cup \{0\} \quad \mbox{and}
  \quad {\bf Z}_0^-:= {\bf Z}\setminus {\bf N}.$$    
For convenience in presentation we also introduce suitable normalizing factors into the definition in (\ref{pPsiq})
and define the normalized Wright function
\begin{equation}\label{pPsiqn}
 {}_p{\hat{\bf \Psi}}_q \left[ \begin{array}{r}(\alpha_1,\,A_1),\,\ldots,\,(\alpha_p,\,A_p);\\
                      (\beta_1,\,B_1),\,\ldots,\,(\beta_q,\,B_q); \end{array}z \right] 
=\sum_{k=0}^\infty \,\frac{\prod_{j=1}^p\,\Gamma( \alpha_j + A_j\,k)/\Gamma(\alpha_j)}{\prod_{j=1}^q\, \Gamma ( \beta_j + B_j\,k)/\Gamma(\beta_j)}\, \frac{z^k}{k!}.
    \end{equation}
    
A special case is
\[{}_p{\hat\Psi}_q \left[ \begin{array}{r}(\alpha_1,\,1),\,\ldots,\,(\alpha_p,\,1);\\
                      (\beta_1,\,1),\,\ldots,\,(\beta_q,\,1); \end{array}z \right]  
                  ={}_pF_q \left[ \begin{array}{r} \alpha_1,\,\ldots,\,\alpha_p\,;\\
                       \beta_1,\,\ldots,\,\beta_q \,; \end{array}
                  \,\, z \right]\]
\begin{equation}\label{pFq}
=\sum_{n=0}^\infty \, {(\alpha_1)_n \cdots (\alpha_p)_n
               \over (\beta_1)_n \cdots (\beta_q)_n} \frac{z^n}{ n!},
\end{equation}
where $(a)_k$ is the Pochhammer symbol or rising factorial defined  by $(a)_n=\g(a+n)/\g(a)=a(a+1)\ldots (a+n-1)$.
Here  an empty product is interpreted as $1$. We assume that the variable
$z$, the numerator parameters
 $\alpha_1,$ $\ldots,$ $\alpha_p,$ and the denominator   parameters
 $\beta_1,$ $\ldots,$ $\beta_q$ take on complex values, provided that no zeros appear in the
 denominator of (\ref{pFq}), that is
\[
  \beta_j \in {\bf C} \setminus {\bf Z}_0^-; \,\, j=1,\ldots, q.
\]
For more details of ${}_pF_q$ including its convergence, its various special and limiting cases,
and its further diverse generalizations,
see, for example, \cite{Bailey, Rainville, Srivastava1}.

The beta function $B(\alpha,\beta)$ is defined by
\begin{equation}\label{Beta-ft}
B(x,\, y) = \left\{ \begin{array}{ll} \int_0^1 \, t^{x -1} (1-t)^{y -1} \, dt & (\Re(x)>0; \,\, \Re(y)>0) \\
\\
\dfrac{\Gamma (x) \, \Gamma (y)}{\Gamma (x+ y)} &
(x,\, y \in {\bf C}\setminus {\bf Z}_0^- ). \end{array} \right.
\end{equation}
Let us recall the well-known integral representation of the ${}_2F_1$ function given by (see, e.g., \cite[p. 85]{Rainville} and \cite[p.~65]{Srivastava1})
\begin{equation}\label{eq(8)}
{}_2F_1 \left(a,\, b;\, c; \, z \right) = \frac{1}{B(b,c-b)}\int_{0}^{1} t^{b-1}(1-t)^{c-b-1} (1-zt)^{-a} \,dt
\end{equation}
\[
(\Re(c)>\Re(b)>0,\,\, \arg(1-z)<\pi).
\]

A number of extensions of some familiar special functions have been recently investigated
 (see, for example,  
 \cite{Chaudhry1, Chaudhry2, Chaudhry3, Chaudhry4, Chaudhry5, Chaudhry6, Chaudhry7, Chaudhry8, Ozergin, Shadab}).
 Chaudhry {\it et al.\/} \cite{Chaudhry1} presented  an extension of the beta function as follows
\begin{eqnarray}\label{eq(9)}
B(x,y;p)= \int_{0}^{1}t^{x-1}(1-t)^{y-1}\exp{\left(-\frac{p}{t(1-t)}\right)}\,dt
\end{eqnarray}
\[
 \left(\Re(p)>0;\,\, ~p=0, ~\Re(x)>0, ~\Re(y)>0\right)
\]
and showed that this extension has certain connections with the Macdonald, error and Whittaker functions.
In addition, the following integral representation of the extended hypergeometric function has been introduced in \cite[Eq. (3.2)]{Chaudhry2}
\begin{eqnarray}\label{eq(10)}
F_{p}(a,b;c;z)=\frac{1}{B(b,c-b)}\int_{0}^{1}t^{b-1}(1-t)^{c-b-1}(1-zt)^{-a}\exp{\left(-\frac{p}{t(1-t)}\right)}\,dt
\end{eqnarray}
\[
 \left(p \in {\bf R}^+;\,~p=0, \,\, |\arg~(1-z)|<\pi,\,\,\Re(c)>\Re(b)>0\right).
\]

The Mittag-Leffler function is defined by \cite[p.~261]{DLMF}
\begin{eqnarray}\label{eq(11)}
E_{\alpha}(z) = \sum_{n=0}^{\infty} \frac{z^n}{\Gamma{(\alpha n + 1)}} \quad \left(\alpha \in {\bf R}^+_0,\,\, z \in {\bf C} \right);
\end{eqnarray}
from which it is obvious that
\[E_{0}(z)=\frac{1}{1-z} \,\,\,(|z|<1),\quad E_{1}(z)=e^z \quad \mbox{and} \quad E_2(z)=\cosh z^{1/2}.  \]
Recently, explicit evaluations of Euler-type integrals of the form
\begin{eqnarray}\label{eq(14)}
\int_{0}^{1}u^{a}(1-u)^{b-1}f(u)du
\end{eqnarray}
have been obtained by Ismail and Pitman \cite{Ismail} for some particular functions $f$, especially in the symmetric case $a=b$.
The evaluation of such integrals is relevant to many reduction formulas for hypergeometric functions that generalize the evaluation of some symmetric Euler-type integrals implicit by the result contained in \cite{pitman}.
Motivated by the work of Choi {\it et al.\/} \cite{Choi1} and Ismail and Pitman \cite{Ismail}, we evaluate some 
Euler-type integrals associated with the Mittag-Leffler function $E_{\lambda}[p\xi(t)]$ of the form
\begin{eqnarray}\label{eq(15)}
\int_{a}^{b}(t-a)^{\alpha-1}(b-t)^{\beta-1}[\chi(t)]^{\gamma}~E_{\lambda}[p\xi(t)]\,dt
\end{eqnarray}
for complex parameter $p$ and some particular choices of the functions $\chi(t)$ and $\xi(t)$.
\vspace{0.6cm}

\begin{center}
{\bf 2. \  Some evaluations of Euler-type integrals}
\end{center}
\setcounter{section}{2}
\setcounter{equation}{0}
\renewcommand{\theequation}{\arabic{section}.\arabic{equation}}
In this section, we derive some theorems on the evaluation of the Euler-type integrals 
\begin{eqnarray}\label{eq(16)}
I^{\alpha,\beta,\gamma}_{a,b,\lambda}[\chi(t),\xi(t);p]=\frac{1}{B(\alpha,\beta)}\int_{a}^{b}(t-a)^{\alpha-1}(b-t)^{\beta-1}[\chi(t)]^{\gamma}~E_{\lambda}[p\xi(t)]\,dt
\end{eqnarray}
for some specific functions $\chi(t)$ and $\xi(t)$. In all cases considered $p\in{\bf C}$ is a complex variable.
\newtheorem{theorem}{Theorem}
\begin{theorem}$\!\!\!.$\ For $a=0$,~$b=\gamma=1$,~$\xi(t)=t(1-t)$~and~$\chi(t)=(1-x_{1}t)^{-\alpha_{1}}(1-x_{2}t)^{-\alpha_{2}}$, we have the following integral
\begin{eqnarray}\label{eq(17)}
&&I^{\alpha,\beta,1}_{0,1,\lambda}\left[(1-x_{1}t)^{-\alpha_{1}}(1-x_{2}t)^{-\alpha_{2}},t(1-t);p\right]\nonumber\\
 && =\sum_{m,n=0}^{\infty}\frac{(\alpha)_{m+n}(\alpha_{1})_{m}(\alpha_{2})_{n}x_{1}^{m}x_{2}^{n}}{(\alpha+\beta)_{m+n}m!n!} \ {}_3{\hat \Psi}_{2}\left[\begin{array}{r} (\alpha\!+\!m\!+\!n,1), (\beta,1), (1,1);\\
~\\
(\alpha\!+\!\beta\!+\!m\!+\!n,2), (1,\lambda) ;\end{array}\!\! p \right]\nonumber\\
\end{eqnarray}
for $\Re(\alpha),~\Re(\beta)>0$, ~$\lambda\geq0$ and ~$|x_{1}|,~|x_{2}|<1$, where we recall that ${}_p{\hat\Psi}_q$ is the normalized Wright function defined in (\ref{pPsiqn}). 
\end{theorem}

\noindent {\bf Proof.}\ We have \cite[p. 962(7)]{Ismail}
\begin{equation}\label{eq(18)}
\int_{0}^{1}t^{\alpha-1}(1-t)^{\beta-1}(1-x_{1}t)^{-\alpha_{1}}(1-x_{2}t)^{-\alpha_{2}}dt=B(\alpha,\beta)F_{1}(\alpha;\alpha_{1},\alpha_{2};\alpha+\beta;x_{1},x_{2})
\end{equation}
\[(\Re(\alpha),~\Re(\beta)>0;~\lambda\geq0; ~|x_{1}|,~|x_{2}|<1),\]
where $F_{1}$ is the Appell double hypergeometric series defined by \cite[p.~413]{DLMF}
\[F_1(\alpha;\beta, \beta';\gamma; x,y)=\sum_{m,n=0}^\infty \frac{(\alpha)_{m+n} (\beta)_m (\beta')_n}{(\gamma)_{m+n} m! n!}\,x^my^n\qquad (|x|,\,|y|<1).\]
Then we set $a=0$,~$b=\gamma=1$,~$\xi(t)=t(1-t)$, $\chi(t)=(1-x_{1}t)^{-\alpha_{1}}(1-x_{2}t)^{-\alpha_{2}}$  in (\ref{eq(16)}) and use the series expansion for $E_\lambda[pt(1-t)]$ in (\ref{eq(11)}). After an interchange of the order of integration and summation and use of the above integral followed by simplification using the properties of the gamma function, we obtain (\ref{eq(17)}).

We note that ${}_3{\hat\Psi}_2(0)=1$ when $p=0$ and the right-hand side of (\ref{eq(17)}) reduces to $$F_1(\alpha;\alpha_{1},\alpha_{2};\alpha+\beta;x_{1},x_{2})$$ in accordance with (\ref{eq(18)}).

\begin{theorem}$\!\!\!.$\ For $a=0$,~$b=\gamma=1$,~$\xi(t)=t(1-t)$ and~$\chi(t)=(1-x_{1}t)^{-\alpha_{1}}(1-x_{2}(1-t))^{-\alpha_{2}}$, we have the following integral
\begin{eqnarray}\label{eq(19)}
&&I^{\alpha,\beta,1}_{0,1,\lambda}\left[(1-x_{1}t)^{-\alpha_{1}}(1-x_{2}(1-t))^{-\alpha_{2}},t(1-t);p\right]\nonumber\\
&&  \hskip 20mm  =\sum_{m,n=0}^{\infty}\frac{(\alpha)_m (\beta)_n (\alpha_{1})_{m}(\alpha_{2})_{n}x_{1}^{m}x_{2}^{n}}{(\alpha+\beta)_{m+n}m!n!}\,{}_3{\hat \Psi}_{2}\left[\begin{array}{r} (\alpha\!+\!m,1), (\beta\!+n,1), (1,1);\\
~\\
(\alpha\!+\!\beta\!+\!m\!+\!n,2), (1,\lambda) ;\end{array}\!\! p \right]\nonumber\\
\end{eqnarray}
\[(\Re(\alpha),~\Re(\beta)>0;~\lambda\geq0; ~|x_{1}|,~|x_{2}|<1).\]
\end{theorem}

\noindent{\bf Proof.}\ We have \cite[p. 279(17)]{Srivastava2}
\begin{equation}\label{e19a}
\int_{0}^{1}t^{\alpha-1}(1-t)^{\beta-1}(1-x_{1}t)^{-\alpha_{1}}(1-x_{2}(1-t))^{-\alpha_{2}}dt=B(\alpha,\beta)F_{3}(\alpha,\beta,\alpha_{1},\alpha_{2};\alpha+\beta;x_{1},x_{2})
\end{equation}
\[(\Re(\alpha),\Re(\beta)>0;~\lambda\geq0; ~|x_{1}|,~|x_{2}|<1),\]
where $F_{3}$ is the Appell double hypergeometric series defined by \cite[p.~413]{DLMF}
\[F_3(\alpha, \alpha';\beta, \beta';\gamma;x,y)=\sum_{m,n=0}^\infty \frac{(\alpha)_m (\alpha')_n (\beta)_m (\beta')_n}{(\gamma)_{m+n} m! n!}\,x^m y^n\qquad (|x|,\, |y|<1).\]
Then on setting $a=0$,~$b=\gamma=1$,~$\xi(t)=t(1-t)$ and~$\chi(t)=(1-x_{1}t)^{-\alpha_{1}}(1-x_{2}(1-t))^{-\alpha_{2}}$  in (\ref{eq(16)}), using (\ref{eq(11)}) and the above integral, we obtain upon simplification (\ref{eq(19)}).

\begin{theorem}$\!\!\!.$\ For $\chi(t)= ut+v$ and $\xi(t)=(t-a)(b-t)$, we have the following integral
\[I^{\alpha,\beta,\gamma}_{a,b,\lambda}\left[ut+v,(t-a)(b-t);p\right]=\sum_{m=0}^{\infty}\frac{(-\gamma)_{m}(\alpha)_m}{(\alpha+\beta)_m m!}\,
\left(\frac{-u(b-a)}{au+v}\right)^{m}\hspace{4cm}\]
\begin{equation}\label{eq(21)}
\hspace{4.5cm}\times {}_3{\hat\Psi}_{2}\left[\begin{array}{r} (\alpha\!+\!m,1), (\beta,1), (1,1);\\
~\\
(\alpha\!+\!\beta\!+\!m,2), (1,\lambda) ;\end{array}\!\! p \right]
\end{equation}
\[(\Re(\alpha),~\Re(\beta)>0; \lambda\geq0;  \left|\arg\left(\frac{bu+v}{au+v}\right)\right|< \pi; a\neq b).\]
\end{theorem}

\noindent{\bf Proof.}\ From \cite[p. 263]{Prudnikov1}, we have 
\[
\int_{a}^{b}(t-a)^{\alpha-1}(b-t)^{\beta-1}(ut+v)^{\gamma}dt=B(\alpha,\beta)~{}_2F_{1}\left(-\gamma, \alpha;
\alpha+\beta;-\frac{(b-a)u}{au+v}\right)
\]
\[(\Re(\alpha),~\Re(\beta)>0;~\lambda\geq0; ~\left|\arg\left(\frac{bu+v}{au+v}\right)\right|< \pi; ~a\neq b).\]
On setting $\chi(t)= ut+v$, and $\xi(t)=(t-a)(b-t)$ in (\ref{eq(16)}), using (\ref{eq(11)}) and the above integral and simplifying, we obtain (\ref{eq(21)}).

\begin{theorem}$\!\!\!.$\ For $\chi(t)= b-a+\nu(t-a)+\mu(b-t)$,~$\xi(t)=(t-a)(b-t)/\chi^2(t)$ and $\gamma=-(\alpha+\beta)$, we have following integral
\begin{eqnarray}\label{eq(23)}
&&I^{\alpha,\beta,-(\alpha+\beta)}_{a,b,\lambda}\left[b-a+\nu(t-a)+\mu(b-t),\frac{(t-a)(b-t)}{(b-a+\nu(t-a)+\mu(b-t))^{2}}\right]\nonumber\\
&&  =\frac{(\nu+1)^{-\alpha}(\mu+1)^{-\beta}}{(b-a)}~{}_3{\hat\Psi}_{2}\left[\begin{array}{r} (\alpha,1), (\beta,1),(1,1);\\
~\\
(\alpha\!+\!\beta,2), (1,\lambda) ;\end{array}\!\! \frac{p}{(\nu+1)(\mu+1)} \right]\nonumber\\
\end{eqnarray}
\[(\Re(\alpha),~\Re(\beta)>0;~\lambda\geq0; ~\chi(t)\neq 0; ~a\neq b).\]
\end{theorem}
\noindent{\bf Proof.} \ From \cite[p. 261(3.1)]{Srivastava3} we have
\begin{eqnarray}\label{eq(24)}
\int_{a}^{b}\frac{(t-a)^{\alpha-1}(b-t)^{\beta-1}}{(b-a+\nu(t-a)+\mu(b-t))^{\alpha+\beta}}\,dt=B(\alpha,\beta)\,\frac{(\nu+1)^{-\alpha}(\mu+1)^{-\beta}}{b-a}\nonumber\\
\end{eqnarray}
\[(\Re(\alpha),~\Re(\beta)>0;~\lambda\geq0;~b-a+\nu(t-a)+\mu(b-t)\neq 0;~a\neq b).\]
On setting $\chi(t)= b-a+\nu(t-a)+\mu(b-t)$,~$\xi(t)=(t-a)(b-t)/\chi^2(t)$ and~$\gamma=-(\alpha+\beta)$ in equation (\ref{eq(16)}), using (\ref{eq(11)}) and the above integral and simplifying, we obtain (\ref{eq(23)}).
\vspace{0.6cm}

\begin{center}
{\bf 3. \ An integral associated with a generating function }
\end{center}
\setcounter{section}{3}
\setcounter{equation}{0}
\renewcommand{\theequation}{\arabic{section}.\arabic{equation}}
Let us consider a two-variable generating function $G(x,t)$ which can be expanded in a formal power series of $t$ such that
\begin{eqnarray}\label{eq(25)}
G(x,t)=\sum_{n=0}^{\infty}c_{n}g_{n}(x)t^{n},
\end{eqnarray}
where the coefficient set $\{c_{n}\}_{n=0}^{\infty}$ may contain the parameters of the set $\{g_{n}(x)\}_{n=0}^{\infty}$ but is independent of $x$ and $t$.\\

\begin{theorem}$\!\!\!.$\ \ Let the generating function $G(x,tu^{\delta}(1-u)^{\omega})$ be uniformly convergent under the conditions $u\in(0,1),~\delta,\omega\geq0$ and $\delta+\omega>0$, with $G(x,t)$ defined by (\ref{eq(25)}). Then
\begin{eqnarray}\label{eq(26)}
\int_{0}^{1}u^{r-1}(1-u)^{s-r-1}G(x,tu^{\delta}(1-u)^{\omega})\,E_{\lambda}[pu(1-u)]\,du\hspace{4cm}\nonumber\\
\hspace{4cm}=\sum_{n=0}^{\infty}c_{n}g_{n}(x)t^{n}{_3}\Psi_{2}\left[\begin{array}{r} (r\!+\!\delta n,1), (s\!-\!r\!+\!\omega n,1), (1,1);\\
~\\
(s\!+\!\delta n\!+\!\omega n,2), (1,\lambda) ;\end{array}\  p \right]
\end{eqnarray}
for $s>r>0$, $\lambda\geq0$, $p\in {\bf C}$, where we recall that ${}_p\Psi_q$ is the un-normalized Wright function defined in (\ref{pPsiq}).
\end{theorem}

\noindent{\bf Proof.}\  On using the definition of generating function $G(x,t)$ given in (\ref{eq(25)}) in the left-hand side of (\ref{eq(26)}), we obtain
\[\sum_{n=0}^{\infty}c_{n}g_{n}(x)t^{n}\int_{0}^{1}u^{r+\delta n-1}(1-u)^{s-r+\omega n-1}E_{\lambda}[pu(1-u)]\,du\]
\[=\sum_{n=0}^{\infty}c_{n}g_{n}(x)t^{n}\sum_{k=0}^\infty \frac{p^k}{\Gamma(1+\lambda k)}\int_0^1u^{r+\delta n+k-1}(1-u)^{s-r+\omega n+k-1}du.\]
Evaluation of the integral as a beta function then yields the right-hand side of (\ref{eq(26)}).

\newtheorem{corollary}{Corollary}
\begin{corollary}$\!\!\!.$\ From (\ref{eq(26)}), we obtain
\begin{eqnarray}\label{eq(28)}
\int_{0}^{1}u^{r-1}(1-u)^{r-1}G(x,t(u(1-u))^{\omega})E_{\lambda}[pu(1-u)]\,du\hspace{4cm}\nonumber\\
\hspace{4cm}=\sum_{n=0}^{\infty}c_{n}g_{n}(x)t^{n}{_3}\Psi_{2}\left[\begin{array}{r} (r+\omega n,1), (r+\omega n,1), (1,1);\\
~\\
(2r+2\omega n,2), (1,\lambda) ;\end{array}\ p \right].
\end{eqnarray}
\end{corollary}
\textbf{Proof.} On setting $s=2r$ and $\delta=\omega$  in (\ref{eq(26)}), we obtain (\ref{eq(28)}).
\bigskip

We now derive some explicit evaluations of definite integrals in terms of Wright hypergeometric functions by
making use of Theorem 5 and consideration of the following generating functions.
\vspace{0.3cm}

\noindent {\bf Example 3.1}\ \  Let us consider the generating function \cite[p.~44(8)]{Srivastava3}
\[
G(x,t)=(1-xt)^{-a}=\sum_{n=0}^{\infty}(a)_{n}\frac{(xt)^{n}}{n!}={}_{1}F_{0}\left[\begin{array}{lll}~\,a;\\
~-;\end{array}xt\right]\qquad (|xt|<1).\]
Use of this generating function (with $x=1$) and Theorem 5 yields
\begin{eqnarray}\label{eq(30)}
\int_{0}^{1}u^{r-1}(1-u)^{s-r-1}[1-tu^{\delta}(1-u)^{\omega})]^{-a}E_{\lambda}[pu(1-u)]\,du\hspace{4cm}\nonumber\\
\hspace{4cm}=\sum_{n=0}^{\infty}\frac{(a)_{n}t^{n}}{n!}{_3}\Psi_{2}\left[\begin{array}{r} (r+\delta n,1), (s-r+\omega n,1), (1,1);\\
~\\
(s+\delta n+\omega n,2), (1,\lambda) ;\end{array}\ p \right]
\end{eqnarray}
for $s>r>0$, $\lambda\geq0$, $\delta,\omega\geq0$, $\delta+\omega>0$.
\vspace{0.3cm}

\noindent{\bf Example 3.2}\ \  Let us consider the generating function \cite[p.409(2)]{Prudnikov2}
\begin{equation}\label{eq(31)}
G(x,t)=\sum_{n=0}^{\infty}\frac{(a)_{n}}{(b)_{n}}{}_{1}F_{1}\left[\begin{array}{lll}\,a;\\
b+n;\end{array}x\right]\frac{t^{n}}{n!}=\Phi_{2}[a,a;b;x,t]\qquad(|x|, |t|<\infty),
\end{equation}
where $\Phi_{2}$ is Humbert's confluent hypergeometric series in two variables defined by \cite{Srivastava3}
\[\Phi_2[b_1,b_2;c;x,y]=\sum_{m,n=0}^\infty \frac{(b_1)_m (b_2)_n x^my^n}{(c)_{m+n} m! n!} \qquad (|x|, |y|<\infty).\]
Use of this generating function and Theorem 5 yields
\begin{eqnarray}\label{eq(32)}
\int_{0}^{1}u^{r-1}(1-u)^{s-r-1}\Phi_{2}[a,a;b;x,tu^{\delta}(1-u)^{\omega})]\,E_{\lambda}[pu(1-u)]\,du\hspace{4cm}\nonumber\\
\hspace{2cm}=\sum_{n=0}^{\infty}\frac{(a)_{n}}{(b)_{n}}\frac{t^{n}}{n!}{}_{1}F_{1}\left[\begin{array}{lll}\,a;\\
b+n;\end{array}x\right]{_3}\Psi_{2}\left[\begin{array}{r} (r+\delta n,1), (s-r+\omega n,1), (1,1);\\
~\\
(s+\delta n+\omega n,2), (1,\lambda) ;\end{array}\ p \right]
\end{eqnarray}
for $s>r>0$, $\lambda\geq0$, $\delta,\omega\geq0$, $\delta+\omega>0$, $|x|<\infty$.
\vspace{0.3cm}

\noindent{\bf Example 3.3}\ \  Let us consider the generating function \cite[p.276(1)]{Rainville}
\begin{eqnarray}\label{eq(33)}
G(x,t)=(1-2xt+t^{2})^{-a}=\sum_{n=0}^{\infty}C_{n}^{(a)}(x)t^{n}\qquad(|xt|<1),
\end{eqnarray}
where $C_{n}^{(a)}(x)$ denotes the Gegenbauer or ultraspherical polynomial \cite{Rainville}.
Use of this generating function (with $x=1$), Theorem 5 and the fact that $C_n^{(a)}(1)=(2a)_n/n!$ \cite[Eq.~(18.5.9)]{DLMF}, then yields
\begin{eqnarray}\label{eq(34)}
\int_{0}^{1}u^{r-1}(1-u)^{s-r-1}[1-tu^{\delta}(1-u)^{\omega})]^{-2a}E_{\lambda}[pu(1-u)]\,du\hspace{4cm}\nonumber\\
\hspace{4cm}=\sum_{n=0}^{\infty}\frac{(2a)_{n}t^{n}}{n!}{_3}\Psi_{2}\left[\begin{array}{r} (r+\delta n,1), (s-r+\omega n,1), (1,1);\\
~\\
(s+\delta n+\omega n,2), (1,\lambda) ;\end{array}\ p \right]
\end{eqnarray}
for $s>r>0$, $\lambda\geq0$, $\delta,\omega\geq0$, $\delta+\omega>0$.
\vspace{0.6cm}

\begin{center}
{\bf 4. \ Applications}
\end{center}
\setcounter{section}{4}
\setcounter{equation}{0}
\renewcommand{\theequation}{\arabic{section}.\arabic{equation}}
We derive some interesting results as applications of the integrals discussed in Section 2.
\vspace{0.3cm}

\noindent{\bf Example 4.1}\ \  
First we take $\alpha_{1}=\alpha_{2}$, $\alpha=\beta>0$ and $x_{2}=x_{1}/(x_{1}-1)$ in (\ref{eq(17)}) to find
\[I^{\alpha,\alpha,1}_{0,1,\lambda}\left[(1-x_{1}t)^{-\alpha_{1}}\left(1-\frac{x_{1}}{x_{1}-1}t\right)^{-\alpha_{1}},t(1-t);p\right]\]
\begin{equation}\label{eq(34)}   =\sum_{m,n=0}^{\infty}\frac{(\alpha)_{m+n}(\alpha_{1})_{m}(\alpha_{1})_{n}x_{1}^{m}(\frac{x_{1}}{x_{1}-1})^{n}}{(2\alpha)_{m+n} m!n!}\, {_3}{\hat \Psi}_{2}\left[\begin{array}{r} (\alpha,1), (\alpha+m+n,1), (1,1);\\
~\\
(2\alpha+m+n,2), (1,\lambda) ;\end{array}\ p \right]
\end{equation}
for $\Re(\alpha)>0$, $\lambda\geq0$, $|x_{1}|<1$, $\Re(x_1)<\fs$.

If $\alpha_{2}=0$ in (\ref{eq(17)}), we obtain
\begin{equation}\label{eq(36)}
I^{\alpha,\beta,1}_{0,1,\lambda}\left[(1-x_{1}t)^{-\alpha_{1}},t(1-t); p\right]
=\sum_{m=0}^{\infty}\frac{(\alpha)_m(\alpha_{1})_{m}x_{1}^{m}}{(\alpha+\beta)_m m!} \,{}_3{\hat\Psi}_{2}\left[\begin{array}{r} (\alpha+m,1), (\beta,1), (1,1);\\
~\\
(\alpha+\beta+m,2), (1,\lambda) ;\end{array}\ p \right].\nonumber\\
\end{equation}
If, in addition, $\alpha_{1}=0$ we find
\begin{equation}\label{eq(37)}
I^{\alpha,\beta,1}_{0,1,\lambda}[1, t(1-t); p]=
{}_3{\hat\Psi}_{2}\left[\begin{array}{r} (\alpha,1), (\beta,1), (1,1);\\
~\\
(\alpha+\beta,2), (1,\lambda) ;\end{array}\  p \right],\nonumber\\
\end{equation}
which reduces to the form
\begin{eqnarray}\label{eq(38)}
I^{\alpha,\beta,1}_{0,1,1}[1, t(1-t); p]&=&{}_{2}F_{2}\left[\begin{array}{c}\alpha,\beta;\\
\frac{\alpha+\beta}{2}, \frac{\alpha+\beta+1}{2};\end{array}\frac{p}{4}\right]
\end{eqnarray}
when $\lambda=1$. The conditions applying to (\ref{eq(36)}) -- (\ref{eq(38)}) are those appearing in Theorem 1
.

\vspace{0.3cm}

\noindent{\bf Example 4.2}\ \  If $x_2=x_1/(x_1-1)$ in (\ref{eq(19)}), we obtain
\[I^{\alpha,\beta,1}_{0,1,\lambda}\left[\frac{(1-x_{1}t)^{-\alpha_{1}-\alpha_2}}{1-x_1}, t(1-t); p\right]
=\frac{(1-x_1)^{-1}}{B(\alpha,\beta)}\int_0^1 t^{\alpha-1}(1-t)^{\beta-1} (1-x_1t)^{-\alpha_1-\alpha_2}\,E_\lambda[pt(1-t)]\,dt.\]
Then replacing $\alpha_1$ by $\alpha_1+\alpha_2$ and $\alpha_2$ by $0$ in the right-hand side of (\ref{eq(19)}), we obtain
\[I^{\alpha,\beta,1}_{0,1,\lambda}\left[\frac{(1-x_{1}t)^{-\alpha_{1}-\alpha_2}}{1-x_1}, t(1-t); p\right]\hspace{8cm}\] 
\begin{equation}\label{eq(39)}
=\frac{1}{1-x_1}\sum_{m=0}^{\infty}\frac{(\alpha)_m(\alpha_{1}+\alpha_2)_m x_{1}^{m}}{(\alpha+\beta)_m m!}\,{}_3{\hat\Psi}_{2}\left[\begin{array}{r} (\alpha+m,1), (\beta,1), (1,1);\\
~\\
(\alpha+\beta+m,2), (1,\lambda) ;\end{array}\  p \right]
\end{equation}
for $\Re(\alpha)>0, \Re(\beta)>0$, $\lambda\geq0$, ~$|x_{1}|<1$. Alternatively, this result can be obtained directly from the above integral by series expansion of $E_\lambda[pt(1-t)]$ and use of the integral representation of the Gauss hypergeometric function in (\ref{eq(8)}).

If $\lambda=1$, (\ref{eq(39)}) reduces to
\[I^{\alpha,\beta,1}_{0,1,1}\left[\frac{(1-x_{1}t)^{-\alpha_{1}-\alpha_2}}{1-x_1}, t(1-t); p\right]\hspace{8cm}\] 
\begin{eqnarray}\label{eq(40)}
=\frac{1}{1-x_1}\sum_{m=0}^{\infty}\frac{(\alpha)_m(\alpha_{1}+\alpha_2)_m x_{1}^{m}}{(\alpha+\beta)_m m!}\,
{}_2F_2\left[\begin{array}{c}\alpha+m, \beta;\\ \frac{\alpha+\beta+m}{2}, \frac{\alpha+\beta+m+1}{2};\end{array}\frac{p}{4}\right]
\end{eqnarray}
under the same conditions as (\ref{eq(39)}).
\vspace{0.3cm}

\noindent{\bf Example 4.3}\ \  
If $\gamma=-\alpha_{1}$, $u=-x_{1}$, $v=1$, $a=0$ and $b=1$ in (\ref{eq(21)}), we obtain
\[I^{\alpha,\beta,-\alpha_{1}}_{0,1,\lambda}\left[(1-x_{1}t),t(1-t); p\right]\hspace{4cm}\]
\begin{equation}\label{eq(41)}  \hspace{4cm} =\sum_{m=0}^{\infty}\frac{(\alpha)_m(\alpha_{1})_{m}(-x_1)^m}{(\alpha+\beta)_m m!} \,{}_3{\hat\Psi}_{2}\left[\begin{array}{r} (\alpha+m,1), (\beta,1), (1,1);\\
~\\
(\alpha+\beta+m,2), (1,\lambda) ;\end{array}\ p \right]
\end{equation}
for $\Re(\alpha),~\Re(\beta)>0$, $\lambda\geq0$, $|x_{1}|<1$.

\vspace{0.3cm}

\noindent{\bf Example 4.4}\ \  
If $\mu=\nu=0$ in (\ref{eq(23)}), we obtain
\[I^{\alpha,\beta,-(\alpha+\beta)}_{a,b,\lambda}\left[b-a,\frac{(t-a)(b-t)}{(b-a)^2}; p\right]\hspace{4cm}\]
\begin{equation}\label{eq(43)}
 \hspace{4cm} =\frac{1}{(b-a)}~{}_3{\hat\Psi}_{2}\left[\begin{array}{r} (\alpha,1), (\beta,1), (1,1);\\
~\\
(\alpha+\beta,2), (1,\lambda) ;\end{array}\  p \right]
\end{equation}
for $\Re(\alpha),~\Re(\beta)>0$, $\lambda\geq0$, $a\neq b$.

\vspace{0.3cm}

\noindent{\bf Example 4.5}\ \  
Finally, if we take $\alpha=\beta$, $a=0$ and $b=1$ in (\ref{eq(23)}), we obtain
\[I^{\alpha,\alpha,-2\alpha}_{0,1,\lambda}\left[1+\nu t+\mu(1-t),\frac{t(1-t)}{(1+\nu t+\mu(1-t))^{2}}; p\right]\hspace{4cm}\]
\begin{equation}\label{eq(44)}
\hspace{4cm}=[(\nu+1)(\mu+1)]^{-\alpha}~{}_3{\hat\Psi}_{2}\left[\begin{array}{r} (\alpha,1), (\alpha,1), (1,1);\\
~\\
(2\alpha,2), (1,\lambda) ;\end{array}\  \frac{p}{(\nu+1)(\mu+1)} \right]
\end{equation}
for $\Re(\alpha)>0$, $\lambda\geq0$.

For $\lambda=1$, (\ref{eq(44)}) reduces to
\[I^{\alpha,\alpha,-2\alpha}_{0,1,1}\left[1+\nu t+\mu(1-t),\frac{t(1-t)}{(1+\nu t+\mu(1-t))^{2}}\right]\hspace{4cm}\]
\begin{equation}\label{eq(45)}
 \hskip 20mm =[(\nu+1)(\mu+1)]^{-\alpha}~{}_{1}F_{1}\left[\begin{array}{r}\alpha;\\
\alpha+\frac{1}{2};\end{array}\frac{p}{4(\nu+1)(\mu+1)}\right]
\end{equation}
for $\Re(\alpha)>0$.
\vspace{0.6cm}

\begin{center}
{\bf 5. \ Concluding remarks}
\end{center}
\setcounter{section}{5}
\setcounter{equation}{0}
\renewcommand{\theequation}{\arabic{section}.\arabic{equation}}
In this paper, we have evaluated some definite Euler-type integrals. We observe that these evaluations are obtained in terms of the Wright hypergeometric function ${}_3\Psi_{2}$. In addition, we have derived a theorem and applied it to obtain evaluations of these integrals related to ${}_3\Psi_{2}$. Further, we remark that the results presented in this paper can be extended to the multi-variable case.

For this purpose, we employ the following integral representation \cite[p. 965(20)]{Ismail} for Lauricella's multiple hypergeometric series in $n$-variables $F_{D}^{(n)}$ \cite{Srivastava3}
\begin{eqnarray}\label{eq(46)}
\int_{0}^{1}t^{\alpha-1}(1-t)^{\beta-1}\prod_{i=1}^{n}(1-x_{i}t)^{-\alpha_{i}}dt=B(\alpha,\beta)F_{D}^{(n)}(\alpha, \alpha_{1},\dots,\alpha_{n};\alpha+\beta;x_{1},\dots,x_{n})
\end{eqnarray}
\[(\Re(\alpha),\Re(\beta)>0, \max\{|x_{1}|,\dots,|x_{n}|\}<1).\]
If we let $\phi(t)=\prod_{i=1}^{n}(1-x_{i}t)^{-\alpha_{i}}$,  $\xi(t)=t(1-t)$, $a=0$ and $b=\gamma=1$ in  (\ref{eq(16)}), and use (\ref{eq(46)}), we obtain upon simplification
\[I^{\alpha,\beta,1}_{0,1,\lambda}\left[\prod_{i=1}^{n}(1-x_{i}t)^{-\alpha_{i}},t(1-t); p\right]=\frac{1}{B(\alpha,\beta)}\sum_{m_{1},\dots,m_{n}=0}^{\infty}\frac{(\alpha_{1})_{m_{1}}\dots(\alpha_{n})_{m_{n}}x_{1}^{m_{1}}\dots x_{n}^{m_{n}}}{m_{1}!\dots m_{n}!}\]
\begin{equation}\label{eq(47)}
\hspace{4cm} \times {_3}\Psi_{2}\left[\begin{array}{r} (\alpha+m_{1}+\dots+m_{n},1), (\beta,1), (1,1);\\
~\\
(\alpha+\beta+m_{1}+\dots+m_{n},2), (1,\lambda) ;\end{array}\ p \right]\nonumber\\
\end{equation}
for $\Re(\alpha),~\Re(\beta)>0$, $\lambda\geq0$, $\max\{|x_{1}|,~\dots,~|x_{n}|\}<1$.
We observe that (\ref{eq(47)}) reduces to (\ref{eq(17)}) for $\alpha_{3}=\alpha_{4}=\dots=\alpha_{n}=0$.

As a consequence, we can extend Theorem 5 as follows:
\begin{theorem}$\!\!\!.$\ \  Let us take the conditions for $G(x,t)$ defined by (\ref{eq(25)})  to be the same as in Theorem 5. Then we have
\[\int_{0}^{1}u^{r-1}(1-u)^{s-r-1}\prod_{i=1}^{j}(1-x_{i}u)^{-\alpha_{i}}G(x,tu^{\delta}(1-u)^{\omega})E_{\lambda}[pu(1-u)]\,du\hspace{4cm}\]
\[
\hspace{2.3cm} =\sum_{n,m_{1},\dots,m_{j}=0}^{\infty}c_{n}g_{n}(x)t^{n}
 \times \frac{(\alpha_{1})_{m_{1}}\dots(\alpha_{n})_{m_{n}}x_{1}^{m_{1}}\dots x_{n}^{m_{n}}}{m_{1}!\dots m_{n}!}\]
 \begin{equation}\label{eq(48)}
\hspace{4cm} \times{_3}\Psi_{2}\left[\begin{array}{r} (r+\delta n+m_{1}+\dots+m_{j},1), (s-r+\omega n,1), (1,1);\\
~\\
(s+\delta n+\omega n+m_{1}+\dots+m_{j},2), (1,\lambda) ;\end{array}\ p \right]
\end{equation}
for $\Re{(\delta)}>\Re{(\omega)}>0$, $\lambda\geq0$, $\delta, \omega>0$, $\delta+\omega>0$, $p\in{\bf C}$.
\end{theorem}
The proof of this theorem is similar to that of Theorem 5 and will be omitted.

\begin{corollary}$\!\!\!.$\ \  By considering Theorem 6 then  we find, when $\lambda=1$,
\[\int_{0}^{1}u^{r-1}(1-u)^{s-r-1}\prod_{i=1}^{j}(1-x_{i}u)^{-\alpha_{i}}G(x,tu^{\delta}(1-u)^{\omega})E_{\lambda}[pu(1-u)]\,du\hspace{4cm}\]
\[\hspace{2.3cm}=\sum_{n,m_{1},\dots,m_{j}=0}^{\infty}c_{n}g_{n}(x)t^{n}
\frac{(\alpha_{1})_{m_{1}}\dots(\alpha_{n})_{m_{n}}x_{1}^{m_{1}}\dots x_{n}^{m_{n}}}{m_{1}!\dots m_{n}!}\]
\begin{equation}\label{eq(49)}
\hspace{3.6cm} \times{_2}\Psi_{1}\left[\begin{array}{r} (r+\delta n+m_{1}+\dots+m_{j},1), (s-r+\omega n,1);\\
~\\
(s+\delta n+\omega n+m_{1}+\dots+m_{j},2);\end{array}\ p \right]
\end{equation}
for 
$\Re{(\delta)}>\Re{(\omega)}>0$; $\delta, \omega>0$; $\delta+\omega>0$; $p\in{\bf C}$.
\end{corollary}

It is clear that by applying Theorem 6 and Corollary 2 to the generating functions $G(x,t)$ defined by (\ref{eq(31)}) and (\ref{eq(33)}), together with some other generating functions, we may obtain many interesting evaluations.

\end{document}